# Security Constrained Unit Commitment with Corrective Transmission Switching

Arun Venkatesh Ramesh
*Student Member, IEEE*
Department of Electrical and Computer Engineering
University of Houston
Houston, TX, USA
aramesh4@uh.edu

Xingpeng Li
*Member, IEEE*
Department of Electrical and Computer Engineering
University of Houston
Houston, TX, USA
xli82@uh.edu

*Abstract*—Traditionally, power system operations use a static network to deliver power and meet demand optimally. Network topology reconfiguration through transmission switching (TS) has gained significant interest recently to reduce the operational cost of power system operations. However, implementation of TS also causes large disturbance in the network and as a result the use of corrective transmission switching (CTS) in response to power system contingencies is currently being researched extensively. This paper emphasizes the importance of CTS to accomplish flexible transmission in *N-1* security-constrained unit commitment (SCUC) model. An *N-1* SCUC mathematical model implementing a dynamic network in the post-contingency scenario is proposed as opposed to current industry practices of static network in short-term operations. The proposed model is tested and validated on the IEEE 24-bus system. The proposed model results in cost-effective implementation and leads to overall reduced cost, and congestion reduction in the post-contingency scenario.

## NOMENCLATURE

| | |
|---|---|
| $g$ | Generator index. |
| $k$ | Transmission element (line or transformer) index. |
| $t$ | Time period index. |
| $n$ | Bus index. |
| $c$ | Line contingency index. |
| $C$ | Set of non-radial transmission contingencies. |
| $G$ | Set of generators. |
| $T$ | Set of Time intervals. |
| $g(n)$ | Set of generators connected bus $n$. |
| $\delta^+(n)$ | Set of lines with bus $n$ as receiving bus. |
| $\delta^-(n)$ | Set of lines with bus $n$ as sending bus. |
| $UT_g$ | Minimum up time for generator $g$. |
| $DT_g$ | Minimum down time for generator $g$. |
| $c_g$ | Linear cost for generator $g$. |
| $c_g^{NL}$ | No-load cost for generator $g$. |
| $c_g^{SU}$ | Start-up cost for generator $g$. |
| $P_g^{min}$ | Minimum capacity of generator $g$. |
| $P_g^{max}$ | Maximum capacity of generator $g$. |
| $R_g^{hr}$ | Regular hourly ramping limit of generator $g$. |
| $R_g^{SU}$ | Start-up ramping limit of generator $g$. |
| $R_g^{SD}$ | Shut-down ramping limit of generator $g$. |
| $R_g^{10}$ | 10-minute outage ramping limit of generator $g$. |
| $P_k^{max}$ | Long-term thermal line limit for line $k$. |
| $b_k$ | Susceptance of line $k$. |
| $P_k^{emax}$ | Emergency thermal line limit for line $k$. |
| $M$ | Real number with huge value |
| $P_{g,t}$ | Output of generator $g$ in time period $t$. |
| $u_{g,t}$ | Commitment status of generator $g$ in time period $t$. |
| $v_{g,t}$ | Start-up variable of generator $g$ in time period $t$. |
| $r_{g,t}$ | Reserve from generator $g$ in time period $t$. |
| $P_{k,t}$ | Line flow of line $k$ in time period $t$. |
| $\theta_{n,t}$ | Phase angle of bus $n$ in time period $t$. |
| $\theta_{m,t}$ | Phase angle of bus $m$ in time period $t$. |
| $d_{n,t}$ | Predicted demand of bus $n$ in time period $t$. |
| $P_{g,c,t}$ | Output of generator $g$ in time period $t$ after outage of line $c$. |
| $P_{k,c,t}$ | Line flow of line $k$ in time period $t$ after outage of line $c$. |
| $\theta_{n,c,t}$ | Phase angle of bus $n$ in time period $t$ after outage of line $c$. |
| $\theta_{n,c,t}$ | Phase angle of bus $n$ in time period $t$ after outage of line $c$. |
| $z_{c,t}^k$ | Line status variable of line $k$ after outage of line $c$ in time period $t$. |
| $Z_{max}$ | The maximum number of transmission elements that are allowed to switch off in each period. |

## I. INTRODUCTION

Several national level directives in recent years stress on the importance to develop a smarter electrical grid. With the influx of renewable technologies and distributed generation, smarter algorithms are required to utilize the flexible network efficiently and reliably. This includes the development of transmission technologies for optimizing the use of transmission [1]. One such venue of introducing the concept of transmission switching (TS) is in SCUC.

Independent system operators (ISOs) collect bids from generators and utilities each day and solve the Day-ahead (DA) SCUC to provide the optimal commitment of day-ahead schedules for generators to meet the predicted load for each hour of the day. A reliable grid is maintained through adherence to system security constraints. The security criteria include, but not limited to, line thermal limits, generator physical limits, reserve requirements and ramping constraints. The grid must also operate reliably in case of emergencies such as transmission line outage or generator loss and an *N*-1 reliable system stresses the importance of modelling such post-contingency situations.

Traditionally, the reliability of the grid is taken care by committing extra generators to ensure ramping and reserve requirements. When transmission long-term thermal limit constraints are active, it leads to congestion in the network. This would force buying power from expensive generating units thereby increasing the cost of operations. However, transmission assets are treated as static networks in this scope and switching of power system elements can help ISOs to relieve network congestion, maintain the system security, and reduce operation costs.

The reliable operation of system network at optimal cost is of prime concern. Therefore, congestion in transmission network needs to be addressed through TS and the research in [1] and [2] points to cost-saving as a result. It can be noted that the infrastructure to perform transmission switching (TS) already exists and this makes it easier to implement without additional investments. However, the use of transmission as a controllable asset today is limited and left to the operator to decide and relieve the network congestion in emergency situations. Typically, voltage profiles and transfer capability of the network are maintained by transmission congestion management schemes in the case of a congested network.

PJM utilizes protocol to let operators remove key lines that are not required for reliability constraints, [3], whereas ISO New England uses the protocol of removing internal transmission line, [4]. Currently, there are no support for operators such as decision support tools to implement network reconfiguration and the current ISO model does not include switching of transmission lines during short-term operations since transmission assets are treated as a static network. This greatly reduces the potential of an efficient implementation of this technology.

TS poses increased problem complexity and most published work looks at various ways to introduce the concept in the mathematical formulation. Significant benefits are noted as a result of utilizing optimal transmission switching problem by co-optimizing the generation and the network topology. The use of TS can improve the market surplus as well as to relieve congestions but also cause significant large system disturbance, [2].

Extensive research has been carried out and have been published to include network reconfiguration for transmission flexibility in SCUC or DC optimal power flow (DCOPF). Mainly TS as a corrective action shows promising results. The use of TS for control of power system is provided by [5] where search techniques were used to solve the problem. However, the concept of corrective transmission switching (CTS) is implemented in [6]. Prior research in [6]-[7] shows that CTS is effective to reduce transmission line overloading.

The complexity introduced due to the additional constraints of network reconfiguration requires better algorithms to provide feasible solutions. A fast CTS for post-contingency was proposed by [8], where generators and network reconfigurations are co-optimized. However, this method only considers limited switching configurations rather than optimal configurations. The use of CTS has been studied in both day-ahead scenario and real time scenarios. In day-ahead scenario, [9] proposes breaking the SCUC algorithm in to unit commitment master problem and TS sub-problem to solve large scale power system network iteratively. In real-time scenario, reference [10] proposed the use of CTS in AC power flow model for real-time contingency analysis by generating candidate switching solutions for each contingency and finally obtaining the top switching solutions for post-contingency violation reduction. Further in [11], the use of CTS in $N$-1-1 SCUC violations through sequential approach has been modelled. It can also be noted from the above literature that use of CTS in $N$-1 SCUC is not completely studied.

Hence in this paper we propose the use of corrective CTS in post-contingency situations in $N$-1 SCUC formulation. While the problem is complex, implementing CTS in day-ahead operations provides sufficient time to obtain feasible solutions. This can be a viable option to mitigate or eliminate the transmission flow violations during contingent scenarios. Since TS may lead to a large system disturbance, the use of TS is limited as a corrective action in the event of an emergency such as transmission line outage. Prior research, [2], also states the cost-saving benefits of co-optimized methods. Therefore, a co-optimized $N$-1 SCUC with optimal TS in post-contingency scenario is proposed as the model.

The rest of this paper is organized as follows. Section II depicts the process of CTS. Section III presents the proposed model of day-ahead $N$-1 SCUC with CTS. Section IV describes the data used for testing. Results and analysis are discussed in Section V. Section VI is dedicated to conclusions drawn from the results. Finally, Section VII discusses the future scope.

## II. CTS AND INDUSTRIAL EXAMPLE

The proposed action of CTS is described illustratively in Fig. 1 which denotes pre-contingency scenario, Fig. 2 which denotes the post-contingency scenario, and Fig. 3 that shows the corresponding CTS action. Fig. 1 represents a pre-contingency scenario where all line flows are below the line thermal limit. This would be the ideal scenario. Power System network flows and generator output are affected due to contingencies such as loss of generators or transmission lines which result in post-contingency scenario. This would require utilize the reserve margin available from committed generators to ramp up in 10 minutes to meet the demand. Fig. 2 represents a post-contingency scenario when line 3 outage occurs. Hypothetically, this leads to an overload of flow on line 4. Practically, this means that line 4 is congested and is the bottleneck in this scenario to deliver power to the loads. To meet such contingencies more generators are required to be committed. A plausible solution for CTS is to remove line 1 or line 2 from the network which is depicted by Fig. 3 which in turn reduces the flow on line 4. CTS can avoid such problems by utilizing the capacity of available generation.

Some industrial examples of switching and control procedures detailed in [12] were utilized by PJM. PJM used one such method as a corrective response to damage caused during Superstorm Sandy. During the storm, PJM lost 82 bulk electric facilities and the system demand was low which led to overvoltage issues in high voltage line. Several 500 kV transmission lines were opened by PJM as part of TS implementation to mitigate the overvoltage [13].

## III. MATHEMATICAL MODEL

An *N*-1 SCUC formulation ensures the loss of any transmission element outages or generator loss. For the purpose of demonstrating CTS, only transmission element outages are considered. The objective is to minimize the total cost of day-ahead operation subject to physical constraints of the generators and transmission elements while meeting post-contingency constraints. The widely used DC power flow model is considered in the proposed approach.

*Objective:*

$$Min \sum_g \sum_t \left(c_g P_{g,t} + c_g^{NL} u_{g,t} + c_g^{SU} v_{g,t}\right) \quad (1)$$

s.t.:

*Base case modeling of generation:*

$$P_g^{min} u_{g,t} \leq P_{g,t}, \forall g, t \quad (2)$$
$$P_{g,t} + r_{g,t} \leq P_g^{max} u_{g,t}, \forall g, t \quad (3)$$
$$0 \leq r_{g,t} \leq R_g^{10} u_{g,t}, \forall g, t \quad (4)$$
$$\sum_{q \in G} r_{q,t} \geq P_{g,t} + r_{g,t}, \forall g, t \quad (5)$$
$$P_{g,t} - P_{g,t-1} \leq R_g^{hr} u_{g,t-1} + R_g^{SU} v_{g,t}, \forall g, t \quad (6)$$
$$P_{g,t-1} - P_{g,t} \leq R_g^{hr} u_{g,t} + R_g^{SD}(v_{g,t} - u_{g,t} + u_{g,t-1}), \forall g, t \quad (7)$$
$$\sum_{q=t-UT_g+1}^{t} v_{g,q} \leq u_{g,t}, \forall g, t \geq UT_g \quad (8)$$
$$\sum_{q=t+1}^{t+DT_g} v_{g,q} \leq 1 - u_{g,t}, \forall g, t \leq T - DT_g \quad (9)$$
$$v_{g,t} \geq u_{g,t} - u_{g,t-1}, \forall g, t \quad (10)$$
$$0 \leq v_{g,t} \leq 1, \forall g, t \quad (11)$$
$$u_{g,t} \in \{0,1\}, \forall g, t \quad (12)$$

*Base case modeling of power flow:*

$$P_{k,t} - b_k(\theta_{n,t} - \theta_{m,t}) = 0, \forall k, t \quad (13)$$
$$-P_k^{max} \leq P_{k,t} \leq P_k^{max}, \forall k, t \quad (14)$$
$$\sum_{g \in g(n)} P_{g,t} + \sum_{k \in \delta^+(n)} P_{k,t} - \sum_{k \in \delta^-(n)} P_{k,t} = d_{n,t}, \forall n, t \quad (15)$$

*Post-contingency 10-minute ramping restriction on generation and modeling of contingencies:*

$$P_{g,t} - P_{g,c,t} \leq R_g^{10} u_{g,t}, \forall g, c \in C, t \quad (16)$$
$$P_{g,c,t} - P_{g,t} \leq R_g^{10} u_{g,t}, \forall g, c \in C, t \quad (17)$$
$$P_g^{min} u_{g,t} \leq P_{g,c,t}, \forall g, c \in C, t \quad (18)$$
$$P_{g,c,t} \leq P_g^{max} u_{g,t}, \forall g, c \in C, t \quad (19)$$

*Post-contingency modeling of power flow for non-radial lines:*

$$P_{k,c,t} - b_k(\theta_{n,c,t} - \theta_{m,c,t}) + (1 - z_{c,t}^k)M \geq 0, \forall k, c \in C, t \quad (20)$$
$$P_{k,c,t} - b_k(\theta_{n,c,t} - \theta_{m,c,t}) - (1 - z_{c,t}^k)M \leq 0, \forall k, c \in C, t \quad (21)$$
$$-P_k^{emax} z_{c,t}^k \leq P_{k,c,t} \leq z_{c,t}^k P_k^{emax}, \forall k, c \in C, t \quad (22)$$
$$\sum_k (1 - z_{c,t}^k) \leq Z_{max}, \forall k, c \in C, t, \quad Z_{max} \in \{0,1,2..\} \quad (23)$$
$$\sum_{g \in g(n)} P_{g,c,t} + \sum_{k \in \delta^+(n)} P_{k,c,t} - \sum_{k \in \delta^-(n)} P_{k,c,t} = d_{n,t}, \forall n, c \in C, t \quad (24)$$

This mixed integer optimization problem reduces operational cost, (1), subject to base-case generation constraints (2)-(12), base-case power flow constraints (13)-(15), post-contingency generator model (16)-(19), and post-contingency non-radial transmission element model (20)-(24). (2) and (3) represents the generator output min-max limits, (4) and (5) are the reserve requirements, (6) and (7) are the hourly ramping consideration, (8) and (9) are the min-up and min-down time of generators. (10) and (11) shows the start-up variable definition. The generator commitment indication variable are bound by binary integrality constraints as shown in (12). The base-case physical power flow constraints represented by (13) which models the power flow with DC line flow equations, (14) which depicts the long-term line thermal limits and (15) which represents nodal balance.

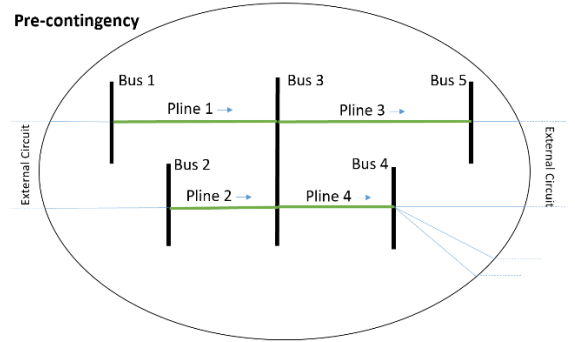

Fig. 1. Line flows in the pre-contingency scenario

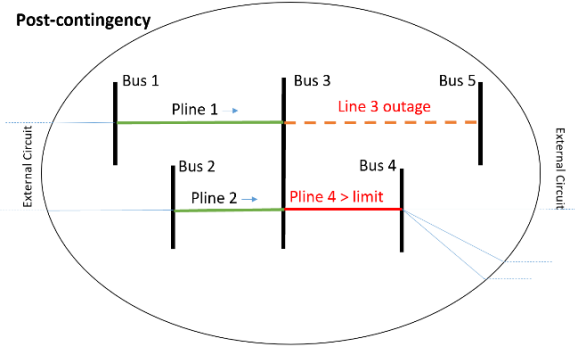

Fig. 2. Line flows in the post-contingency scenario

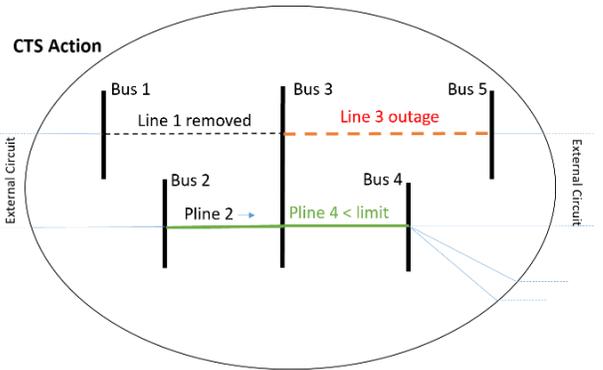

Fig. 3. Post-switching line flows

The mixed integer program is a co-optimization of base-case and post-contingency scenarios. The post-contingency generator constraints are modelled for the base-case solution through (16)-(17) which takes consideration of the 10-minute generator ramping constraints, and (18)-(19), the post-contingency generator output min-max limits when line *c* is lost. Post-contingency scenarios of line flows are modelled for all non-radial lines when a transmission element outage occurs. Equation (22) shows the emergency line flow limits and (24) is the post-contingency nodal balance. The CTS action is represented by the binary decision variable, $z_{c,t}^k$, introduced in equations (20)-(21). This represents the status of switchable transmission element *k* under contingency *c* in time period *t* (When the value is 1, the line is in service and when value 0 the line is switched off/ not is service). *M*, is often represented as 'big M' which is a large real value. It ensures that equation (20) and (21) are linear in nature. The decision variable for switching decides the optimal network configurations for each contingency to relieve the post- contingency congestion in the system. Since TS can cause large system disturbance, restriction on number of transmission elements open is represented in (23).

## IV. TEST CASE: IEEE 24-BUS (RTS 96) SYSTEM

The IEEE 24-bus network developed by power experts [14] was used for testing in this paper. However, a modified data from [15] for the same network was utilized. Fig. 4 represents this modified network which contains 24 buses, 33 generators, and 38 branches. The total generation capacity is 3,393 MW and the types of generator available with operational cost, min and max outputs are presented in Table I. The goal of the proposed SCUC model is to find out cost-saving in congested networks. The system peak load is 2,265 MW with a maximum nodal load of 210 MW and a minimum nodal load of 0.

TABLE I. GENERATOR DATA IN IEEE-24 BUS

| No. of Gen | Min | Max | $/MWh |
|---|---|---|---|
| 4 | 2.4 | 12 | 94.74 |
| 4 | 15.8 | 20 | 163.02 |
| 4 | 15.2 | 76 | 19.64 |
| 6 | 0 | 50 | 0 |
| 3 | 25 | 100 | 75.64 |
| 4 | 54.25 | 155 | 15.46 |
| 3 | 68.95 | 197 | 74.75 |
| 1 | 140 | 350 | 15.89 |
| 2 | 100 | 400 | 5.46 |

## V. RESULTS AND ANALYSIS

The mathematical model is implemented using AMPL and solved using Gurobi solver with a MIPGAP of 0.01 for a 24-hour (Day-Ahead) load period, [16]-[17]. The difference in overall cost of *N*-1 SCUC with CTS and *N*-1 SCUC without CTS is used to demonstrate the cost reduction with CTS.

TABLE II. OPERATIONAL COST IN *N*-1 SCUC

| | *N*-1 SCUC without CTS | | *N*-1 SCUC with CTS | |
|---|---|---|---|---|
| | Scenario I | Scenario II | Scenario I | Scenario II |
| Cost ($) | 932,911 | 921,812 | 923,995 | 921,812 |
| $CC$ ($) | 11,099 | N/A | 2,183 | N/A |

We have modelled the *N*-1 SCUC in two different scenarios. Scenario I is when regular emergency rating and Scenario II is when infinite emergency rating is used for transmission elements in the network respectively. In Scenario I, the operational cost of $932,911 for *N*-1 SCUC without CTS and $923,995 for *N*-1 SCUC with CTS was obtained. This information is tabulated in Table II. Scenario II operational cost shows the operational cost of the system when there are no congestion in the system in post-contingency scenario. This implies that post-contingency congestion is significantly reduced with the use of CTS. It can be further verified from the dual values of the active post-contingency emergency thermal limits constraints; the use of CTS leads lower dual values which is tightly correlated to the congestion cost.

$$CC = TC_{Scenario\ I} - TC_{Scenario\ II} \qquad (25)$$

We define the post-contingency congestion cost ($CC$) in (25). It is defined as the difference in total operational cost in Scenario I ($TC_{Scenario\ I}$) and total operational cost in Scenario II ($TC_{Scenario\ II}$). It can be inferred that in *N*-1 SCUC without CTS, the post contingency congestion cost is $11,099. Whereas in *N*-1 SCUC with CTS, the post-contingency congestion cost is $2,183. For the test scenario, implementing CTS results in reduction of congestion cost by 80.33%. Transmission networks are built with redundancy and by including the CTS the flexibility in the network is utilized which reduces the congestion in the network.

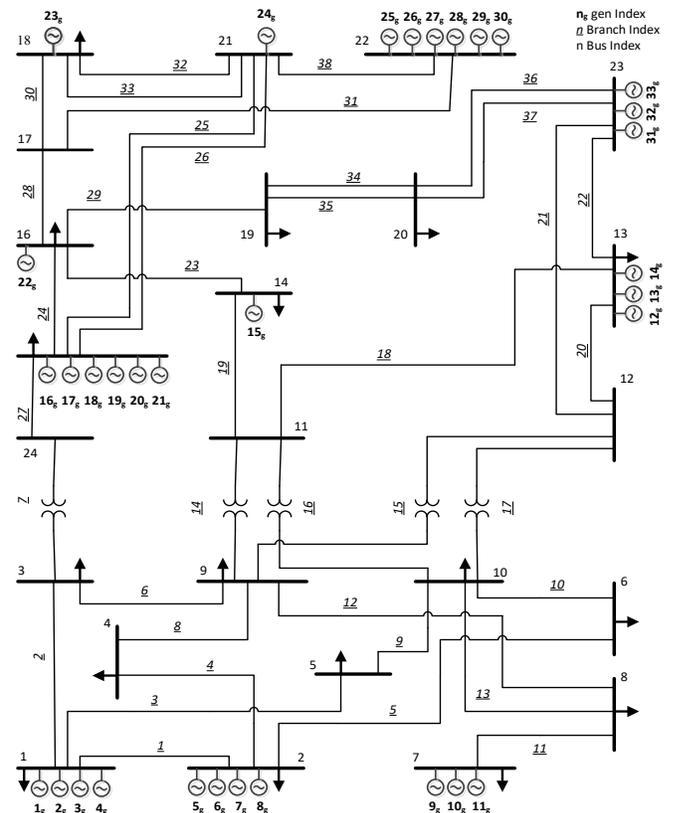

Fig. 4. IEEE 24-bus System – one area of IEEE RTS-96 system [18]

TABLE III. GENERATOR START-UP

| Time Period | ON Generators in *N*-1 SCUC without CTS | ON Generators *N*-1 SCUC with CTS |
|---|---|---|
| t=1 | 3,4,7,8,9,11,21-33 | 3,4,7,8,11,13,21-33 |
| t = 9 | 14 | 9 |
| t = 21 | 16,17,20 | N/A |
| t=23 | 1,5,6 | N/A |

Table III shows that implementing CTS required less frequent generator start-ups in the 24-hour period as TS provides the opportunities for committed generators to ramp-up and meet post-contingency demand. In total, there were 26 generators start-ups in the 24-hour period while using *N*-1 SCUC without CTS whereas only 20 generators start-ups in *N*-1 SCUC with CTS. In both cases 19 generators were started in the first hour. Implementing CTS required only one additional generator start-up as it enabled the existing generators to ramp up without transmission violations to meet the demand.

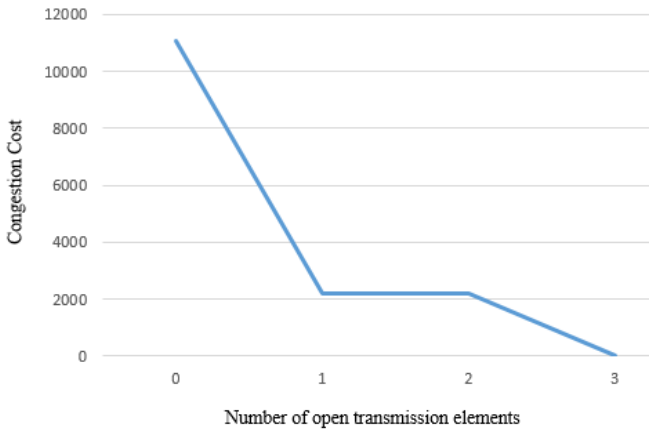

Fig. 5. Number of open transmission elements vs congestion cost

The number of transmission elements opened as part of CTS played a role in reducing the congestion cost. Fig. 5 shows that the congestion cost reduced to $0 from $11,099 when more transmission elements were allowed to be opened if required. The difference in congestion cost when one transmission element is allowed to be opened in CTS versus multiple transmission elements is $2,183.

TABLE IV. Post-contingency congestion scenario for Period 9

| Post-contingency congested line (line number [from-bus – to-bus]) | Post-contingency line outage (line number [from-bus – to-bus]) | |
|---|---|---|
| | *N*-1 SCUC without CTS | *N*-1 SCUC with CTS |
| 10 [6-10] | 1 [1-2],2 [1-3],7 [3-24],8 [4-9],9 [5-10],27 [5-24] | 2 [1-3] |
| 23 [14-16] | 7 [3-24], 18 [11-13],21 [12-13],22 [13-23],27 [5-24] | 7 [3-24] |

The solution of switching transmission lines for each contingent line was studied for a 24-hour period. We can understand the switching ideology using a demand of 2,076 MW for period 8 and a demand of 2,265 MW for period 9. In period 9, after the load profile change, it was noted that line 10 and line 23 are susceptible to post-contingency congestion. The scenario leading to congestion is tabulated in table IV. Line 10 connects from bus 6 to bus 10 with a long-term thermal rating of 157.5 MW and an emergency rating of 180 MW and line 23 connects from bus 14 to bus 16 with a long-term thermal rating of 315 MW and an emergency rating of 393.75 MW. During congestion of line 10 and 23 in period 9, we noticed that the scenarios leading to post-contingency congestion were reduced. CTS was beneficial to produce a maximum line overload reduction of 4% and 24% in lines 10 and 23 respectively. In the best case, 24% reduction of line overload brings the line flow below the long-term thermal limit which reduces significant stress on transmission lines.

TABLE V. SWITCHING SOLUTION FOR PERIOD 9

| CTS action for Period 9 in *N*-1 SCUC with CTS | |
|---|---|
| Outage Line | Switched OFF Line |
| 4 | 31 |
| 7 | 34 |
| 8 | 5 |
| 12 | 37 |
| 17 | 31 |
| 21 | 31 |
| 30 | 15 |
| 32 | 23 |
| 34 | 28 |
| 36 | 31 |
| 37 | 35 |
| 38 | 10 |

In the same period 9, the number of contingent scenarios leading to congestion is 6 and 5 for lines 6 and 10 respectively in *N*-1 SCUC without CTS. In *N*-1 SCUC with CTS, only 1 contingent scenario led to congestion in both lines. The contingent scenario for line 23 is when outage of line 7 occurs. Similarly, line 10 is congested only for outage of line 2. The switching pattern for period 9 is represented in Table V.

VI. CONCLUSIONS

The best scenario is represented by infinite transmission capacity in the post-contingency scenario, which serves as a benchmark to measure the performance of the proposed CTS in SCUC. It is observed that CTS can alleviate the network congestion in post-contingency scenarios by rerouting power through the network. The implementation of CTS also led to fewer generator start-ups. This is evident from the results that only 1 generator start-up is required when CTS is used as compared to 7 without CTS after period 1. Overall, this results in reduced operational cost, congestion cost and higher transmission capability in the case of a congested network.

Studying the line flows in contingent scenarios, we note that line overload was reduced with CTS in most contingent scenarios. The use of CTS can lead to the removal of post-contingency transmission congestions if more transmission elements are allowed to open in each contingent scenarios which will result in $0 in congestion cost. However, there are concerns with TS as it can cause a large disturbance to the system. The additional cost due to the restriction of allowing one transmission element to be open in CTS is a tradeoff between system reliability and cost saving. The congestion cost, $2,183, is only 0.2% of the total operation cost and it can be attributed as a reliability cost to avoid system disturbance.

VII. FUTURE WORK

It can be noted that all transmission elements that are not radial lines are simulated in this co-optimized method. This significantly increases the complexity of the solution. In reality, ISO's model a ranked list of watch list lines or high-risk transmission lines. Such lists can significantly increase the computing time of this model. Also, consideration of generator outages can be considered in the model to provide complete consideration to *N*-1 reliability. Breaking the problem as a Master-Slave problem will improve problem scalability.